\documentclass[10pt]{article}
\usepackage{times,amssymb,amsthm,amsopn,amsfonts,mathrsfs,latexsym,graphicx,citesort}
\usepackage{amsmath}
\usepackage{wrapfig}

\newtheorem{theorem}{Theorem}

\newtheorem{definition}[theorem]{Definition}

\newcommand{\RR}{\mathbb{R}}

\newcommand{\eg}{\varepsilon}

\long\def\symbolfootnote[#1]#2{\begingroup%
\def\thefootnote{\fnsymbol{footnote}}\footnote[#1]{#2}\endgroup}

\begin{document}

\title{Time reversal in thermoacoustic tomography  - an error estimate}
\author{Yulia Hristova\\
\begin{small} \it{ Department of Mathematics,
Texas A\& M University, College Station, TX 77845, USA} \end{small} }
\date{}

\maketitle

\begin{abstract}
In thermoacoustic tomography an object is irradiated by a short electromagnetic pulse and the absorbed energy causes a thermoelastic expansion. This expansion leads to a pressure wave propagating through the object. The goal of thermoacoustic tomography is the recovery of the initial pressure inside the object from measurements of the pressure wave made on a surface surrounding the object. The time reversal method can be used for approximating the initial pressure when the sound speed inside the object is variable (non-trapping as well as trapping). This article presents error estimates for the time reversal method in the cases of variable, non-trapping sound speeds. Numerical examples for non-trapping as well as for trapping sound speeds are provided.
\end{abstract}

\section{Introduction}
Thermoacoustic tomography is a hybrid medical imaging technique characterized by high resolution and contrast. A short electromagnetic pulse is used to irradiate the biological object of interest. A part of the electromagnetic energy is absorbed by the tissue, which causes a thermoelastic expansion. This leads to a pressure wave propagating through the object, which is measured by ultrasound transducers located on an \emph{observation surface} $S$, usually surrounding the object. The collected information is used to reconstruct the initial pressure, which is roughly proportional to the absorbed energy \cite{Tam,MXW_review}. The good contrast in the resulting images is due to electromagnetic energy being preferentially absorbed by cancerous cells, while ultrasound provides a submilimeter resolution \cite{XuWang2002,MXW_review}. 

Let $S$ be a smooth observation surface surrounding the object and let $B$ be the domain bounded by $S$. In applications only space dimensions $n=2, 3$ are of interest, but the analysis developed in this paper is carried out for an arbitrary dimension $n$. We assume that the speed $c(x)$ of ultrasound in the tissue is smooth and strictly positive $c(x)>c>0$, and that $c(x)\equiv 1$ for large values of $|x|$. Then the pressure $p(x,t)$ at location $x$ and time $t$ satisfies the wave equation (e.g. \cite{Wang-book,Diebold,Tam,KuKu}): 
\begin{equation}\label{E:TAT}
\left\{\begin{array}{cc}
p_{tt}=c^2(x)\Delta p, &t\geq 0, \quad x\in\RR^n\\
p(x,0)=f(x), &p_t(x,0)=0, \\
p(y,t)=g(y,t) &\mbox{ for }y\in S, \quad t\geq 0.
\end{array}\right.
\end{equation}
Here $p_{tt}$ denotes the second time derivative of $p(x,t)$, $g(y,t)$ is the measured data, i.e. the value of the pressure at time $t$ measured at transducer's location $y\in S$, and $f(x)$ is the initial pressure, which is to be reconstructed. The function $f(x)$ is assumed to have a compact support in $B$. One thus faces the problem of inverting the mapping $R\colon f\rightarrow g$ from the initial pressure to the measured data. For a more detailed description of the mathematical methods used in thermoacoustic tomography, one can refer to \cite{FinchRak07,KuKu,AKK,FinchRak_inbook,PatchScher07,Haltmeier05} and the references therein. 
There are various types of reconstruction procedures for closed observation surfaces, e.g. filtered backprojection formulas, eigenfunction expansion methods and time reversal method. A comparison of the three and a discussion of their advantages and limitations can be found in \cite{HKN}. It was argued there that time reversal is the most versatile and easy to implement among these. Here we consider the errors involved in the time reversal method and provide estimates that justify its validity. It should be noted that, although we are interested in applications to thermoacoustic tomography, the time reversal method for the wave equation has been used in other applications \cite{FinkPrada01,Fink08,BlomPapanicZhao02}. Our results apply in all these cases.

In Section \ref{S:second} a non-trapping condition on the sound speed is recalled and the time reversal method is described. In Section \ref{S:error_est} the main result of this artice is formulated and proved. Numerical examples are provided in Section \ref{S:numerics}. Section \ref{S:conclusion} contains some concluding remarks. This is followed by the acknowledgment section.

Additional discussions and implementations of the time reversal method can be found in \cite{BurgMattHaltPalt,Grun,HKN}.

\section{Non-trapping condition and time reversal method} \label{S:second}
In this section we give a more detailed description of the time reversal method and recall the notion of a non-trapping sound speed.
\subsection{Non-trapping condition}\label{SS:non_trap}
We will be interested in the initial value problem:
\begin{equation}\label{E:wave}
\left\{\begin{array}{cc}
p_{tt}=c^2(x)\Delta p, &t\geq 0, \quad x\in\RR^n\\
p(x,0)=f_1(x), &p_t(x,0)=f_2(x), \\
\end{array}\right.
\end{equation}
where $c(x)>c>0,\; c(x) \in C^\infty(\RR^n)$, and $c(x)-1$, as well as $f_1(x)$ and $f_2(x)$, has compact support. 
Consider the following Hamiltonian system in $2n$ real variables $(x,\xi)$ with Hamiltonian $H=\frac{c^2(x)}{2}|\xi|^2$:
\begin{equation}\label{E:non-trap}
\left\{\begin{array}{l}
x^\prime_t=\frac{\partial H}{\partial \xi}=c^2(x)\xi \\
\xi^\prime_t=-\frac{\partial H}{\partial x}=-\frac 12 \nabla
\left(c^2(x)\right)|\xi|^2 \\
x|_{t=0}=x_0, \xi|_{t=0}=\xi_0.
\end{array}\right.
\end{equation}
The solutions of this system are called \emph{bicharacteristics} and their projections into the $x$-space $\RR^n_x$ are called \emph{rays}.
\begin{definition}\label{D:non-trap}
We say that the  non-trapping condition holds, if all rays (with $\xi_0\neq 0$) tend to infinity when $t \to\infty$ .
\end{definition}
It is known that the singularities of the solution of (\ref{E:wave}) propagate along bicharacteristics (e.g. \cite{Nirenberg,Taylor_book,EgorovShubin}). Therefore, due to the non-trapping condition imposed on $c(x)$, it follows that for any distributions $f_1,\;f_2$ with compact supports, the singularities of the solution move away to infinity as $t\rightarrow \infty$. Then, for a sufficiently large $t$, the solution $p(x,t)$ is infinitely differentiable inside the unit ball $B$. Moreover, the solution decreases inside $B$ and the following estimates hold:
\begin{theorem}(e.g. \cite{Vainberg,EgorovShubin})\label{T:est} Under the conditions imposed on the sound speed $c(x)$ and
$f(x)$ and for any bounded domain $B$, the solution of (\ref{E:wave}) satisfies the estimates
\begin{equation}
\left|\frac {\partial^{k+|m|}p}{\partial_t^k\partial_x^m}\right| \leq
C\eta_k(t)\left(\|f_1\|_{L^2}+\|f_2\|_{L^2}\right),\; x\in B,\; t>T_0,
\end{equation}
for any multi-index (k,m), where $\eta_k(t)=t^{-n+1-k}$ for even n, and
$\eta_k(t)=e^{-\delta t}$ for odd n. Here $\delta$ is a positive constant depending only on $c(x)$, $T_0$ depends on the domain $B$, and the $C$ depends on $B$ and the multi-index $(k,m)$.
\end{theorem}
We will make use of the above theorem in estimating the error of reconstruction of $f(x)$ when using time reversal. Note that in the case of trapping sound speed the local energy of the solution of (\ref{E:wave}) still decreases in any compact domain, but, in general, there is no uniform local energy decay estimate \cite{Ralston69,Pauen}.

\subsection{The time reversal method}  
The Huygens' principle states that in odd dimensions and when the sound speed is constant, for any initial source with bounded support and for any bounded domain, there is a moment in time when the wave leaves the domain (e.g., \cite{CourHilb}). Thus, given an initial pressure $f(x)$ with a bounded support, there is a time $T$ when the wave inside the domain $B$, bounded by the observation surface $S$, vanishes for all $t\geq T$. Then, to reconstruct $f(x)$ we can simply "rewind" the solution, i.e. we can solve the wave equation backwards in time inside $B$ with zero initial conditions at $t=T$ and boundary conditions given by the measured data on $S$. At time $t=0$ the solution of this problem will be equal to $f(x)$. 

In even dimensions or when the sound speed is variable, the Huygens' principle does not hold anymore. Nevertheless, we could still try to rewind the wave in hopes of approximating the initial pressure. This is the main idea of the \emph{time reversal method} and we will now give a precise description of the method in the general case of a variable sound speed and dimension $n$. 

To reconstruct the initial pressure $f(x)$ inside $B$, we try to reverse the time (starting from time $t=T$ and going back to $t=0$)
and solve the following problem:
\begin{equation}\label{E:inv_exact}
\left\{ \begin{array}{ll}
             u_{tt} = c^2(x)\Delta u \qquad &\mbox{in} \; B\times[0,T] \\
             u(x,T) = p(x,T) & \\
             u_t(x,T)= p_t(x,T)& \\
             u|_{S}(x,t)= g(x,t) \; &\mbox{on}\; S\times [0,T], 
         \end{array}\right.
\end{equation}
where $g = p|_{S\times [0,T]}$ is the restriction of $p$ onto $S\times [0,T]$.
At time $t=0$ the solution of this problem equals $f(x)$\footnote{Note that $u(x,t)$ is simply the restriction of $p(x,t)$ onto $B\times[0,T]$.}. Obviously, this reconstruction requires the knowledge of the solution to (\ref{E:TAT}) on the cylinder $S \times [0,T]$ and inside $B$ at time $T$. The values of $p$ on the cylinder can be obtained from the measurements of the transducers placed on $S$. However, the values of the pressure inside $B$ at time $T$ are not known. Nevertheless, due to Theorem \ref{T:est}, we do know that the solution inside $B$ decays with time. Therefore, after some time $T>T_0$ it is reasonable to approximate $p(x,T)$ and $p_t(x,T)$ with zero. This is how we arrive to the approximate reconstruction problem:
\begin{equation}\label{E:inv_approx}
\left\{ \begin{array}{ll}
             v_{tt} = c^2(x)\Delta v & \mbox{in} \; B\times[0,T] \\
             v(x,T) = 0 &\\
             v_t(x,T)= 0 &\\
             v|_{S}(x,t)= g(x,t)\varphi_\eg(t) & \mbox{on} \;S\times [0,T], 
        \end{array}\right.
\end{equation}
Here $\varphi_\eg(t)$ is a smooth cut-off function that equals to 1 in $(-\infty,T-\eg]$ (where $T-\eg>T_0$) and vanishes for $ t \geq T$. As it will be explained below, at time $t=0$ the solution to (\ref{E:inv_approx}) approximates $f(x)$. 
In Section \ref{S:error_est} we estimate the decay of this error with respect to the cut-off time T in the case of a non-trapping sound speed. 

\section{An error estimate}\label{S:error_est}
Consider the error $e(x,t;T)=u(x,t)-v(x,t)$ that results from replacing the true solution of (\ref{E:inv_exact}), $u(x,t)$, with the time reversal solution  $v(x,t)$. Recall that $v(x,t)$ solves equation (\ref{E:inv_approx}). Then, the error satisfies the following equation: 
\begin{equation}\label{E:err}
\left\{ \begin{array}{l}
             e_{tt}(x,t;T) = c^2(x)\Delta e(x,t;T) \qquad \mbox{in} \; B\times[0,T] \\
             e(x,T;T) = p(x,T)\\
             e_t(x,T;T)= p_t(x,T)\\
             e|_S(x,t;T)= g_\eg(x,t)\colon=(1-\varphi_\eg(t))g(x,t) \qquad \mbox{in}\;S\times[0,T]. 
         \end{array}\right.
\end{equation}

\begin{figure}[!ht]
\begin{center}
\scalebox{0.4}{\includegraphics{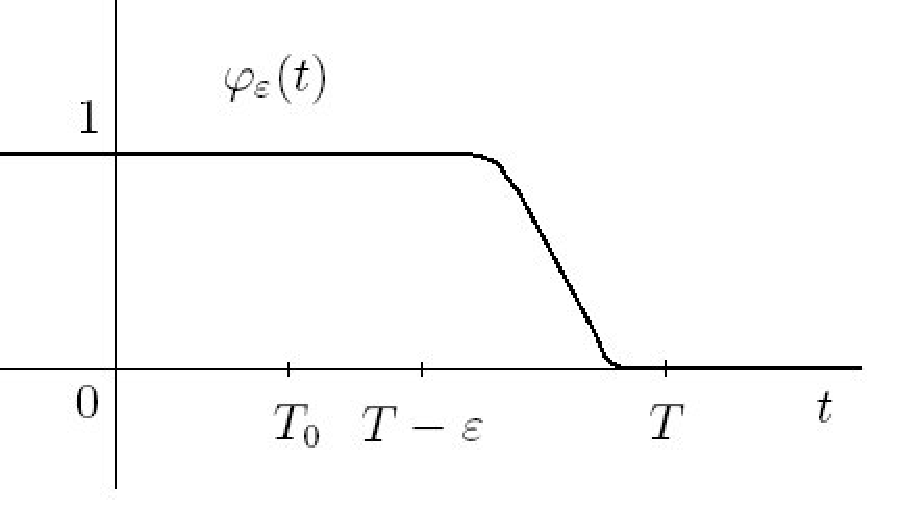}}
\end{center}
\caption{A sketch of a cut-off function $\varphi_\eg(t)$. }
\label{F:phi}
\end{figure}

To make the role of $\eg$ more clear we will specify our choice of a cut-off function $\varphi_\eg(t)$. Let $\varphi(t)$ be a smooth function that equals to $1$ on $(-\infty,-1]$ and vanishes on $[0,\infty)$. When $\eg \leq 1$, we set $\varphi_\eg(t)=\varphi((t-T)/\eg)$ (Figure \ref{F:phi}), and when $\eg>1$, we choose $\varphi_\eg(t)=\varphi_1(t)=\varphi(t-T)$. This can also be written as $\varphi_\eg(t)=\varphi((t-T)/\alpha)$, where $\alpha=\min\{\eg,1\}$. 

Let the sound speed $c(x)$ satisfy the conditions described in Section \ref{SS:non_trap}, namely, $c(x)>c>0,\; c(x) \in C^\infty(\RR^n)$, $c(x)-1$ has a compact support and $c(x)$ is non-trapping. Suppose also that the initial pressure $f(x)$ belongs to $L^2(\RR^n)$ and is compactly supported (in $B$).  Then the following theorem, which is the main result of this article, holds. 
\begin{theorem} \label{T:Theorem}
There exists $T_0$ such that for any $T>T_0$ and $\eg>0$ satisfying $T-\eg>T_0$, the error $e(x,t;T)$ can be estimated as follows:
\begin{itemize}
\item for even dimensions $n$
\begin{equation*}
  \max_{0\leq t\leq T} (\|e(t;T)\|_{H^1(B)}+\|e_t(t;T)\|_{L^2(B)}) \leq  C(\eg)(T-\eg)^{-n+1}\|f\|_{L^2(B)};
\end{equation*}
\item for odd dimensions $n$
\begin{equation*}
  \max_{0\leq t\leq T} (\|e(t;T)\|_{H^1(B)}+\|e_t(t;T)\|_{L^2(B)}) \leq C(\eg) e^{-\delta(T-\eg)}\|f\|_{L^2(B)}.
\end{equation*}
\end{itemize}
Here $C(\eg)=C/\min\{\eg^2,1\}$, for some constant $C$ depending on $B$ and $c(x)$.

In particular, 
\[\|e(0;T)\|_{H_0^1(B)}=\|f-v(0)\|_{H^1(B)} \leq C(\eg)(T-\eg)^{-n+1}\|f\|_{L^2} \mbox{, for even } n;\]
\[\|e(0;T)\|_{H_0^1(B)}=\|f-v(0)\|_{H^1(B)} \leq C(\eg) e^{-\delta(T-\eg)}\|f\|_{L^2(B)} \mbox{, for odd } n,\]
where $v(0,x)$ is the approximation to $f(x)$ obtained by time reversal with cut-off time $T$.  
\end{theorem}
{\bf Proof.~~}
Let $E\colon H^s(S)\rightarrow H^{s+1/2}(B)\;,s>0,$ be the operator of harmonic extension that produces a harmonic function $E\phi$ in $B$ from the Dirichlet boundary data $\phi$ (e.g. \cite{LionsI}). Then $w\colon=e-Eg_\eg$ satisfies
\begin{equation}
\left\{ \begin{array}{l}\label{E:hom_bdry}
             w_{tt} - c^2(x)\Delta w = -{Eg_\eg}_{tt}\qquad \mbox{in} \; B\times[0,T] \\
             w(x,T) = p(x,T)-Eg_\eg(x,T) = p(x,T)-Eg(x,T)\\
             w_t(x,T)= p_t(x,T)-{Eg_\eg} _t(x,T) = p_t(x,T)-Eg_t(x,T)\\
             w|_{S}(x,t) = 0 \quad \mbox{on} \;S\times [0,T], 
         \end{array}\right.
\end{equation}  
As noted in Section \ref{SS:non_trap}, there exists a time $T_0$ after which the solution $p(x,t)$ to (\ref{E:TAT}) is infinitely smooth in $\bar{B}\times [T,\infty)$ for any $T>T_0$. Let us choose such $T$ and let $\eg$ be such that $T-\eg>T_0$. Then the right hand side and the initial conditions of (\ref{E:hom_bdry}) are infinitely smooth functions. Indeed, $p(x,T)$ is smooth, as we discussed above. Recall that $g_\eg$ is a smooth cut-off of the trace of $p(x,t)$. Moreover, it is supported in $[T-\eg,T]$ -- a time interval where $p(x,t)$ has already become infinitely smooth. Thus $Eg_\eg(x,t)$ is a smooth function as well.

Next, we will apply the following theorem to the solution of (\ref{E:hom_bdry}). 
\begin{theorem}(e.g.\cite{Evans})\label{T:Evans}
Assume that $G, H\in C^\infty (\bar{B}), F\in C^\infty(\bar{B}\times[0,T])$ and consider the initial boundary value problem
\begin{equation}\label{E:U}
\left\{\begin{array}{ll}
U_{tt} - c^2(x)\Delta U = F \qquad &\mbox{in} \;  B\times (0,T]\\  
U = G,\quad U_t = H \quad &\mbox{on} \; B \times \{t = 0\},\\
U = 0 \quad &\mbox{on} \;  S\times [0, T]
\end{array}\right.
\end{equation}
Here $c(x)\in C^\infty(\bar{B})$ is such that $c(x) \geq \theta$ for some $\theta>0$.  Let the following compatibility conditions hold for $l=1,2,...$ .
\begin{equation}\label{compat_cond}
\begin{array}{l}
G_0:=G\in H^1_0(B),\;H_1 :=H\in H^1_0(B),\\
G_{2l} :=\frac{d^{2l-2}F}{dt^{2l-2}}(\cdot,0)+c^2(x)\Delta G_{2l-2}\in H^1_0(B),\\
H_{2l+1} :=\frac{d^{2l-1}F}{dt^{2l-1}}(\cdot,0)+c^2(x)\Delta H_{2l-1}\in H^1_0(B).
\end{array}
\end{equation}
Then (\ref{E:U}) has a unique solution $U \in C^\infty(\bar{B}\times[0,T])$ and 
\[\begin{array}{l}
\displaystyle\max_{0\leq t\leq T} (\|U(t)\|_{H_0^1(B)}+\|U_t(t)\|_{L^2(B)})\leq\smallskip\\
 \mbox{\hskip 120pt}Ce^{C_1T}(\|F\|_{L^2(0,T;L^2(B))} +\|G\|_{H_0^1(B)} + \|H\|_{L^2(B)}),    
\end{array}\]
where the constants $C$ and $C_1$ depend on $B$ and $c(x)$.
\end{theorem}
We will verify that the conditions of Theorem \ref{T:Evans} are satisfied by (\ref{E:hom_bdry}), noting that the compatibility conditions must hold at $t=T$, since the direction of time is reversed. Indeed, as we already discussed, $F=-{Eg_\eg}_{tt}$, $G=p(x,T)-Eg(x,T)$ and $H=p_t(x,T)-Eg_t(x,T)$ are infinitely smooth. It is clear that $G_0=G$ and $H_1=H$ belong to $H_0^1(\bar{B})$, as $g(x,t)$ has been defined to be the trace of $p(x,t)$. The higher order compatibility conditions hold as well: 
\[\begin{array}{l@{\,}l}
G_{2l}(x)=&-\frac{d^{2l-2}}{dt^{2l-2}}{Eg_\eg}_{tt}(x,T)+c^2\Delta G_{2l-2}=-E\frac{d^{2l}}{dt^{2l}}{g}(x,T)+c^2\Delta G_{2(l-1)}\smallskip\\
& = -E\frac{d^{2l}}{dt^{2l}}{g}(x,T)+(c^2\Delta)^l G_0 = -E\frac{d^{2l}}{dt^{2l}}{g}(x,T)+(c^2\Delta)^l p(x,T).
\end{array}\]
We used here that $\Delta E = 0$.

As we already know that $p(x,T)$ is smooth in a neighborhood of $T$, we conclude that 
\[\left. G_{2l}(x)\right|_S=\left.\left[-\frac{d^{2l}}{dt^{2l}}{p}(x,T)+(c^2\Delta)^l p(x,T)\right]\right|_S = 0,\]
so $G_{2l}\in H^1_0(B)\mbox{ for }l=1,2...$ Similarly, one can check that $H_{2l+1}\in H^1_0(B)\mbox{ for }l=1,2...$.

Before applying  Theorem \ref{T:Evans} to the solution of (\ref{E:hom_bdry}) we notice that the energy in $B$, defined by \mbox{$\mathscr{E}(t) := \frac{1}{2}\int_B \left( \left|\nabla w\right|^2 +c^{-2}(x)\left|w_t\right|^2\right)dx$}, stays constant in $[0,T-\alpha]$, where $\alpha=\min\{\eg,1\}$. Indeed, using (\ref{E:hom_bdry}) one easily shows that
\[
\dot{\mathscr{E}}(t) = \int_B \left( \nabla w\cdot\nabla w_t + c^{-2}w_t w_{tt}\right)\;dx = -\int_B w_t c^{-2}E{g_{\eg}}_{tt}\;dx.\]

The last integral vanishes in $[0,T-\alpha]$, because ${g_\eg}_{tt}\equiv 0$ in $ [0,T-\alpha]$ and therefore, $E{g_{\eg}}_{tt}\equiv0$ in this interval. Then, by Theorem (\ref{T:Evans}), it follows that
\begin{equation*}\label{E: with energy}
\begin{array}{l@{\!}l}
\displaystyle\max_{0\leq t \leq T}\mathscr{E}(t)= \max_{T-\alpha \leq t \leq T}\mathscr{E}(t)\leq 
  C\max_{T-\alpha \leq t \leq T}&\left\{ \|w(t)\|^2_{H^1_0(B)}+\|w_t(t)\|^2_{L^2(B)} \right\} \leq\\
 C\left\{ \|{Eg_\eg}_{tt}\|^2_{L^2(T-\alpha,T;L^2(B))} \right.
  + \|p(\cdot,T) - &Eg(\cdot,T)\|^2_{H_0^1(B)} \\
  &+ \left. \|p_t(\cdot,T)-Eg_t(\cdot,T)\|^2_{L^2(B)}\right\},
\end{array}
\end{equation*}
where $C$ depends on $c(x)$ and $B$ only, but not on $T$ or $\eg$. Here and in what follows $C$ will denote various constants, all of them independent on $T$ or $\eg$.

As $w=e-Eg_\eg$, the above inequalities imply 
\begin{equation} \label{eq:err}
\begin{array}[t]{l}
  \displaystyle\max_{0\leq t\leq T}(\|e(t)\|_{H^1(B)}+\|e_t(t)\|_{L^2(B)}) \leq\\
 C \displaystyle\max_{0\leq t \leq T}\bigg\{\sqrt{\mathscr{E}(t)}+\|Eg_\eg(t)\|_{H^1(B)} + \|{Eg_\eg}_t(t)\|_{L^2(B)} \bigg\} \leq\\
  \displaystyle C\bigg\{\max_{T-\alpha\leq t\leq T}(\|Eg_\eg(t)\|_{H^1(B)}+\|{Eg_\eg}_t(t)\|_{L^2(B)})+\|{Eg_\eg}_{tt}\|_{L^2(T-\alpha,T;L^2(B))} \\
  \mbox{\hskip 50pt} +\|p(\cdot,T)-Eg(\cdot,T)\|_{H_0^1(B)}  + \|p_t(\cdot,T)-Eg_t(\cdot,T)\|_{L^2(B)}\bigg\}.
\end{array}
\end{equation}
We used here that $Eg_\eg \equiv 0$ in $[0,T-\alpha]$.

Let $\mathcal{T}\colon H^s(B)\rightarrow H^{s-1/2}(S), s>1/2,$ denote the trace operator, which is known to be continuous (e.g.\cite{Evans}). 
Then, using the continuity of the linear operators $E$ and $\mathcal{T}$, we can estimate the terms in (\ref{eq:err}) as follows.
\[ \begin{array}{lr}
\|Eg_\eg(t)\|_{H^1(B)}\leq & C\|{g_\eg(t)}\|_{H^{1/2}(S)} = C\|(1-\varphi_\eg(t))\mathcal{T}p(t)\|_{H^{1/2}(S)}\leq \smallskip\\
&C\|(1-\varphi_\eg(t))p(t)\|_{H^1(B)}\leq C\|p(t)\|_{H^1(B)}.
\end{array}\]
Similarly, 
\[ \begin{array}{lr}
\|{Eg_\eg}_t(t)\|_{L^2(B)}\leq &\|{Eg_\eg}_t(t)\|_{H^1(B)}\leq C\|[(1-\varphi_\eg)p(t)]_t\|_{H^1(B)}\leq \smallskip\\
&\displaystyle\frac{C}{\alpha}\left\{ \|p(t)\|_{H^1(B)}+\|p_t(t)\|_{H^1(B)} \right\},
\end{array}\]
\[ \begin{array}{r}
\|{Eg_\eg}_{tt}(t)\|_{L^2(B)} \leq \|{Eg_\eg}_{tt}(t)\|_{H^1(B)} \leq C\|[(1-\varphi_\eg)p]_{tt}\|_{H^1(B)} \leq \smallskip\\
 \displaystyle\frac{C}{\alpha^2}\left\{ \|p(t)\|_{H^1(B)}+\|p_t(t)\|_{H^1(B)}+\|p_{tt}(t)\|_{H^1(B)}\right\},
\end{array}\]
\[
\|{Eg(T)}\|_{H^1(B)}\leq C\|{p(T)}\|_{H^1(B)} \mbox{ and } \|{Eg_t(T)}\|_{L^2(B)}\leq C\|{p_t(T)}\|_{H^1(B)}.
\]
Here we used that $\displaystyle\max_t |\varphi_\eg^\prime(t)|\leq C^\prime/\alpha $ and $\displaystyle\max_t |\varphi_\eg^{\prime\prime}(t)|\leq C^{\prime\prime}/\alpha^2$ for some constants $C^\prime$ and $ C^{\prime\prime}$.
Then, the error can be estimated by
\begin{equation}\label{E:interm}
\begin{array}{l}
\displaystyle\max_{0\leq t\leq T}\left\{\|e(t)\|_{H^1(B)}+\|e_t(t)\|_{L^2(B)}\right\}\leq \smallskip\\
\displaystyle \frac{C}{\alpha^2}\left\{ \max_{T-\alpha \leq t\leq T} \left(\|p(t)\|_{H^1(B)}+\|p_t(t)\|_{H^1(B)}\right)+\right.\smallskip\\
\mbox{\hskip 15pt}\quad\|p\|_{L^2(T-\alpha,T;H^1(B))}+ \|p_t\|_{L^2(T-\alpha,T;H^1(B))} + \|p_{tt}\|_{L^2(T-\alpha,T;H^1(B))} + \smallskip\\
\mbox{\hskip 200pt}\|p(T)\|_{H^1(B)} +\|p_t(T)\|_{H^1(B)}\Bigg\}
\end{array}
\end{equation}
Theorem \ref{T:est} allows us to estimate the quantity inside the braces in right-hand side of above inequality by a factor of
\[\begin{array}{l}
 \displaystyle \bigg\{\max_{T-\alpha\leq t\leq T}\left[\eta_0(t)+\eta_1(t)\right]+\|\eta_0\|_{L^2(T-\alpha,T)}+ \|\eta_1\|_{L^2(T-\alpha,T)}\smallskip\\
\mbox{\hskip 100pt} + \|\eta_2\|_{L^2(T-\alpha,T)} + \eta_0(T) + \eta_1(T)\bigg\}\|f\|_{L^2(B)}
\end{array}\]
The functions $\eta_k(t)$ are monotonically decreasing, and $\eta_2(t)\leq\eta_1(t)\leq\eta_0(t)$ provided $t\geq1$. 
Therefore, the right hand side of (\ref{E:interm}) is less than 
\[C/\alpha^2\left\{ \eta_0(T-\alpha)+\|\eta_0\|_{L^2(T-\alpha,T)}\right\}\|f\|_{L^2(B)}.\] 
Taking into account that 
\[\|\eta_0\|^2_{L^2(T-\alpha,T)}\leq \alpha\max_{T-\eg\leq t\leq T}\eta^2_0(t) \leq \eta^2_0(T-\eg),\] we arrive to 
\[
\max_{0\leq t\leq T}(\|e(t)\|_{H^1(B)}+\|e_t(t)\|_{L^2(B)})\leq C(\eg)\eta_0(T-\eg)\|f\|_{L^2(B)},
\]
where $\eta_0(T-\eg)=(T-\eg)^{-n+1}$ for even n, $\eta_0(T-\eg)=e^{-\delta(T-\eg)}$ for odd n, and $C(\eg)=C/\min\{\eg^2,1\}$.

This proves the theorem.
\section{Numerics of errors}\label{S:numerics}

\begin{figure}[!ht]
\begin{center}
 \scalebox{0.3}{\includegraphics{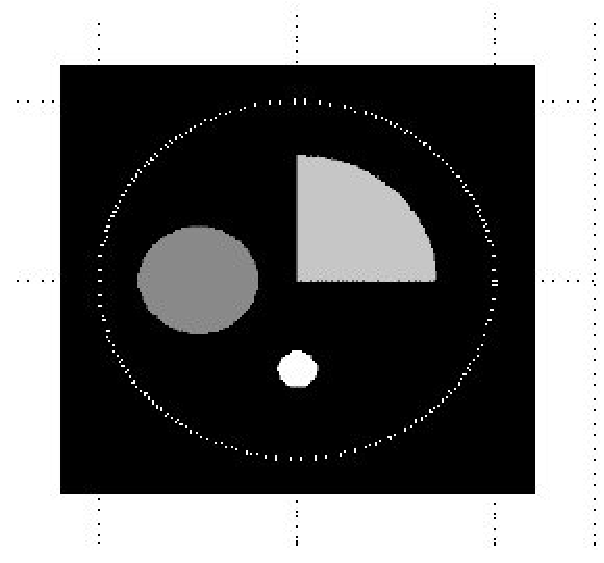}}\vspace{0.2cm}\\
\scalebox{0.27}{\includegraphics{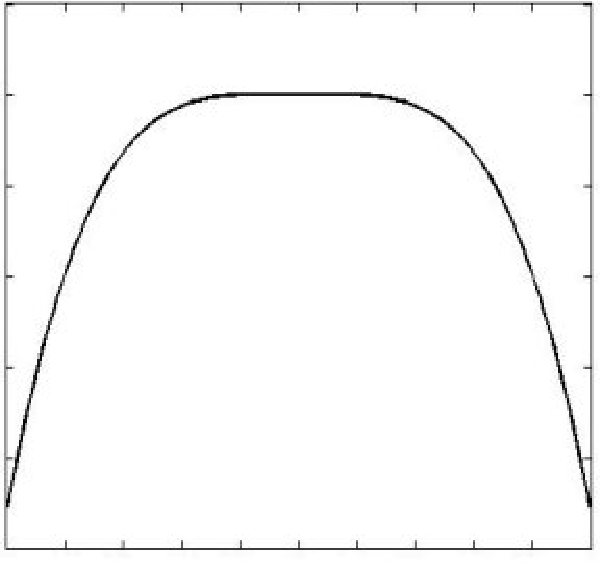}} \mbox{\hskip 50pt} \scalebox{0.3}{\includegraphics{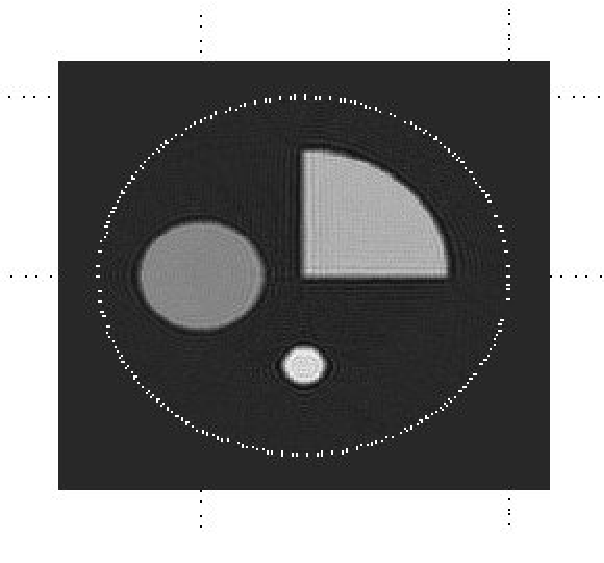}}\vspace{0.2cm}\\
\scalebox{0.27}{\includegraphics{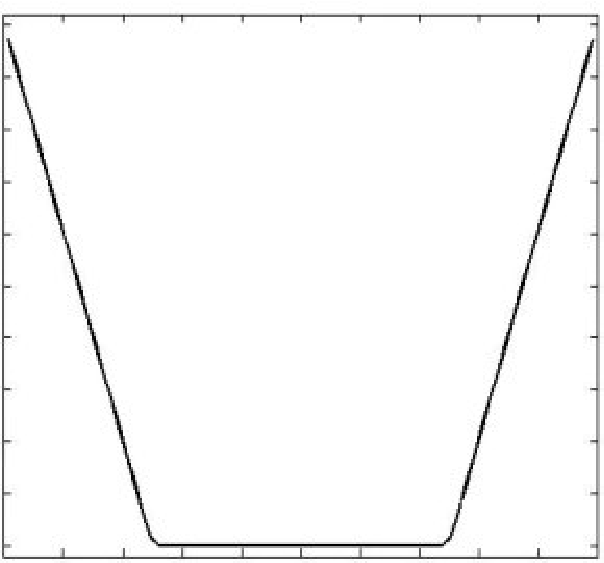}} \mbox{\hskip 50pt} \scalebox{0.3}{ \includegraphics{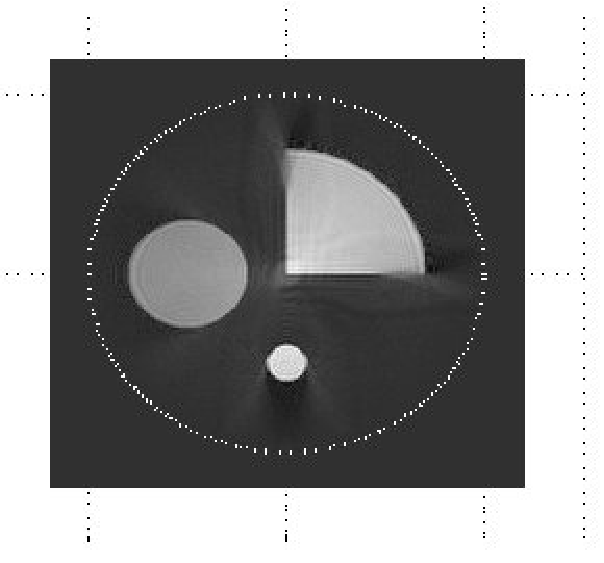}}\vspace{0.2cm}\\
\scalebox{0.27}{\includegraphics{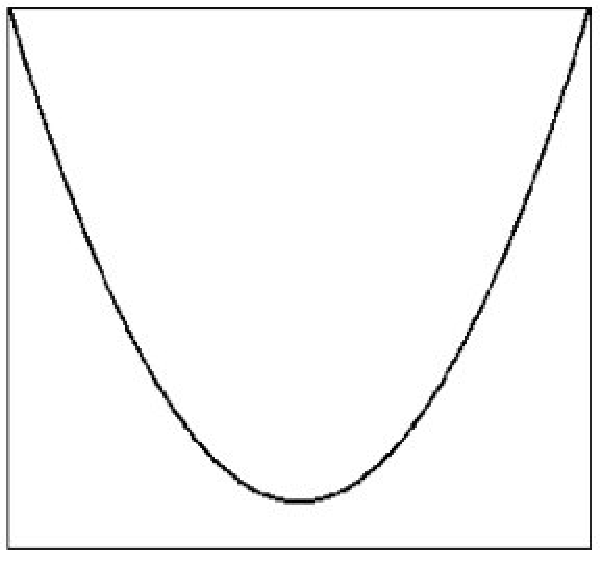}} \mbox{\hskip 50pt} \scalebox{0.3}{ \includegraphics{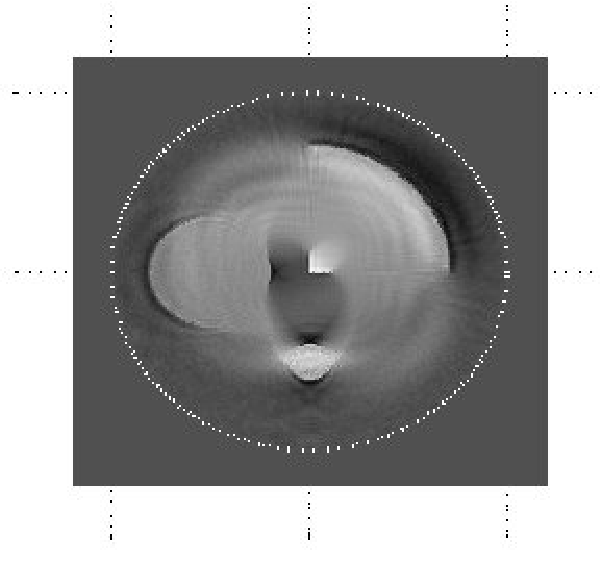}}
\end{center}
\caption{A phantom (top) and its reconstructions using various radial sound speeds. Profiles of radial sound speeds are shown in the left column. In the right column are the reconstructions of the phantom, obtained by using the corresponding sound speeds from the left column. The first sound speed (second row) is non-trapping. The second one is a ``trapping crater`` sound speed. The third speed is a paraboloid, which causes severe trapping. The white dotted circles represent the observation surface $S$. }
\label{F:rec}
\end{figure}

\begin{figure}[!ht]
\begin{center}
 \scalebox{0.3}{\includegraphics{cos_sp.eps}}\vspace{0.2cm}\\
\scalebox{0.25}{\includegraphics{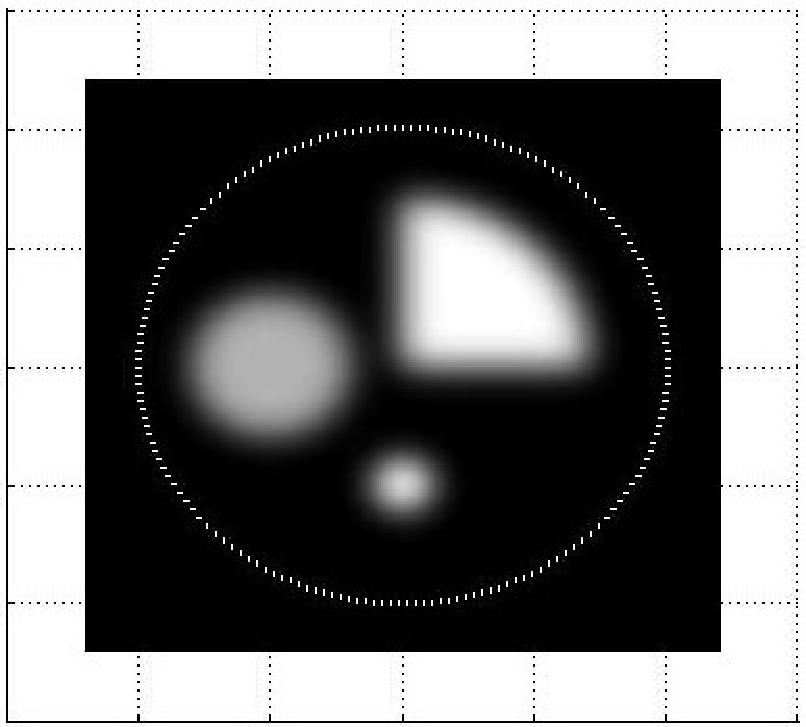}} \mbox{\hskip 50pt} \scalebox{0.27}{\includegraphics{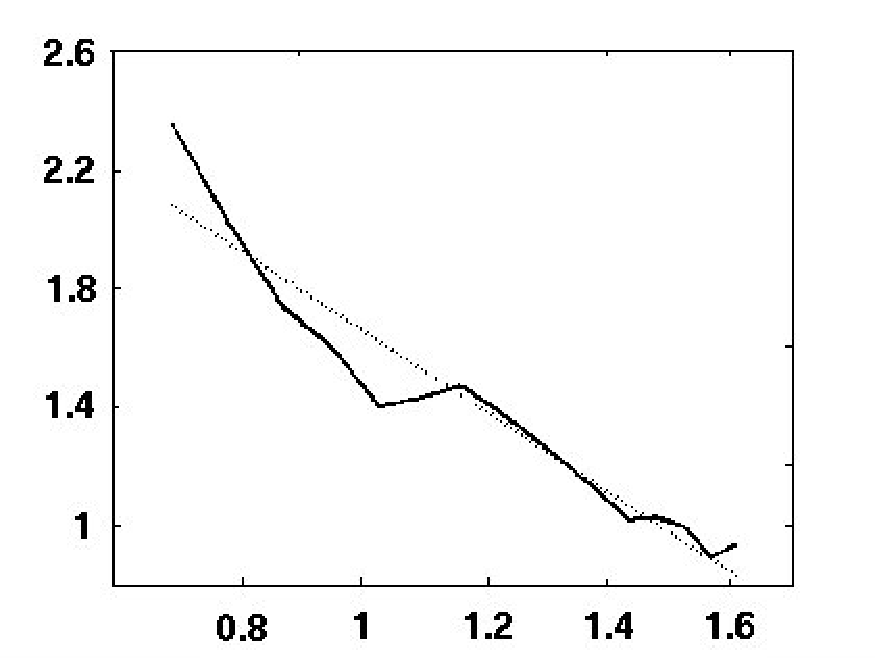}}\vspace{0.2cm}\\
\scalebox{0.25}{\includegraphics{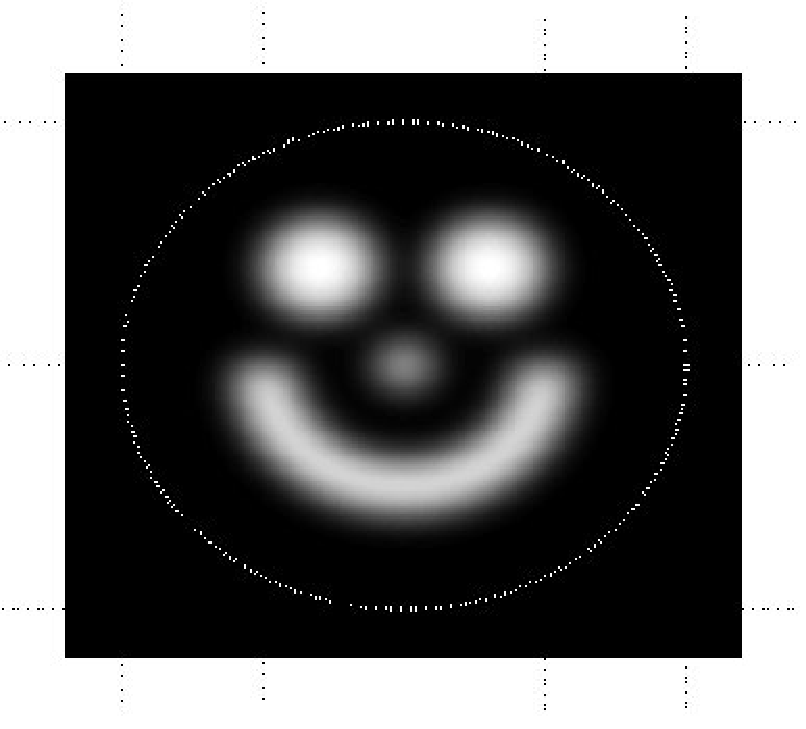}} \mbox{\hskip 50pt} \scalebox{0.27}{ \includegraphics{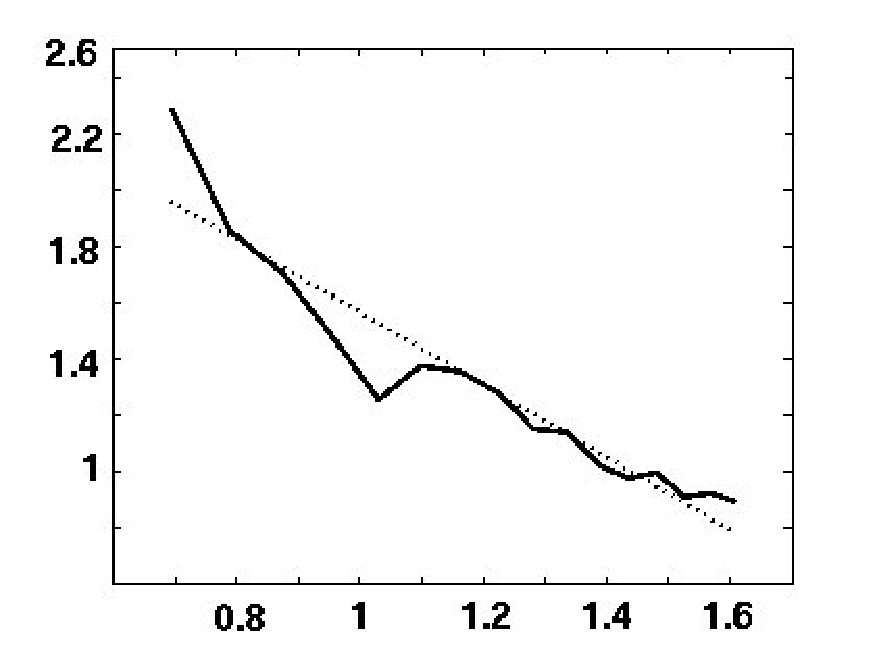}}
\end{center}
\caption{Axial profile of a radial non-trapping sound speed (top), density plot of phantoms (left column) and the corresponding $H^1$ errors of reconstruction as functions of the cut-off time $T$ (right column). The plots of the errors are in logarithmic scale, i.e. the horizontal axis represents $\ln T$ and the vertical axis represents $\ln (H^1$ error$)$. The dotted lines show the linear regression interpolation of the error. The decay of the error in the first example (second row) is of the order of $T^{-1.36}$, in the second one (third row) it is $T^{-1.28}$.}  
\label{F:nice_cos}
\end{figure}
\begin{figure}[!ht]
\begin{center}
\scalebox{0.25}{\includegraphics{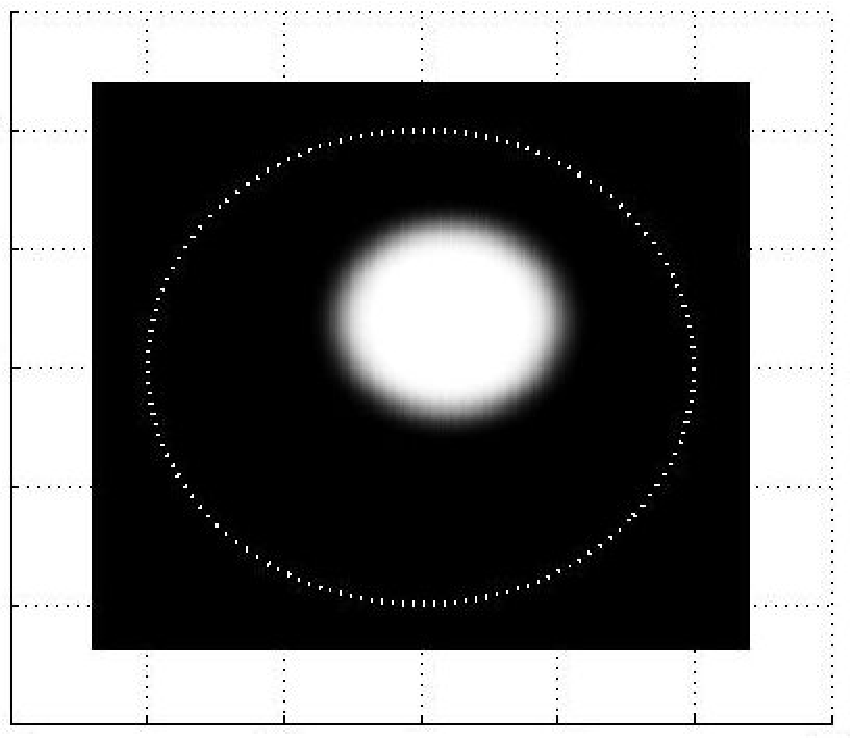}}\vspace{0.2cm}\\
\scalebox{0.25}{\includegraphics{vsm_comb1_ph.eps}} \mbox{\hskip 50pt} \scalebox{0.27}{\includegraphics{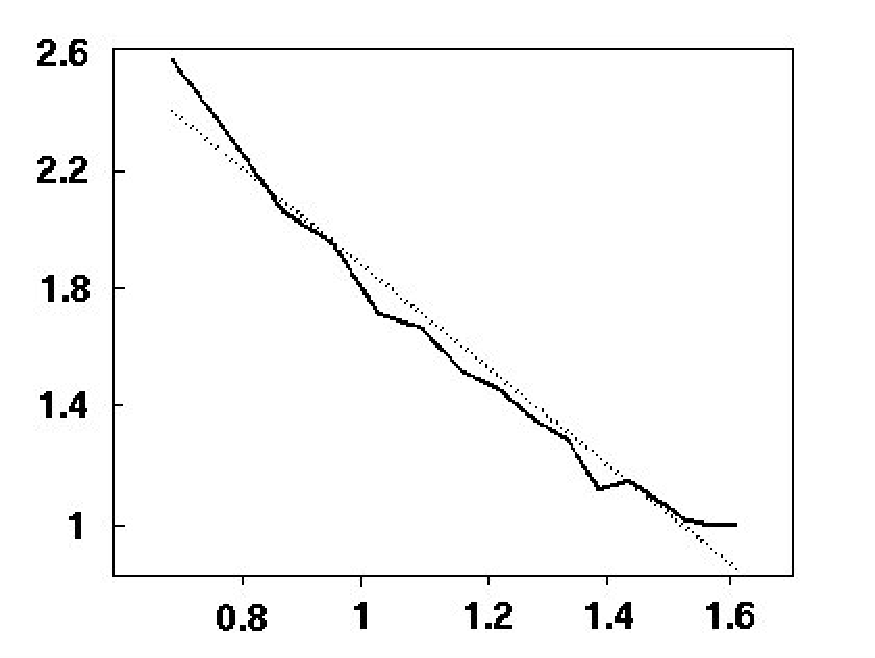}}\vspace{0.2cm}\\
\scalebox{0.25}{\includegraphics{face.eps}} \mbox{\hskip 50pt} \scalebox{0.27}{\includegraphics{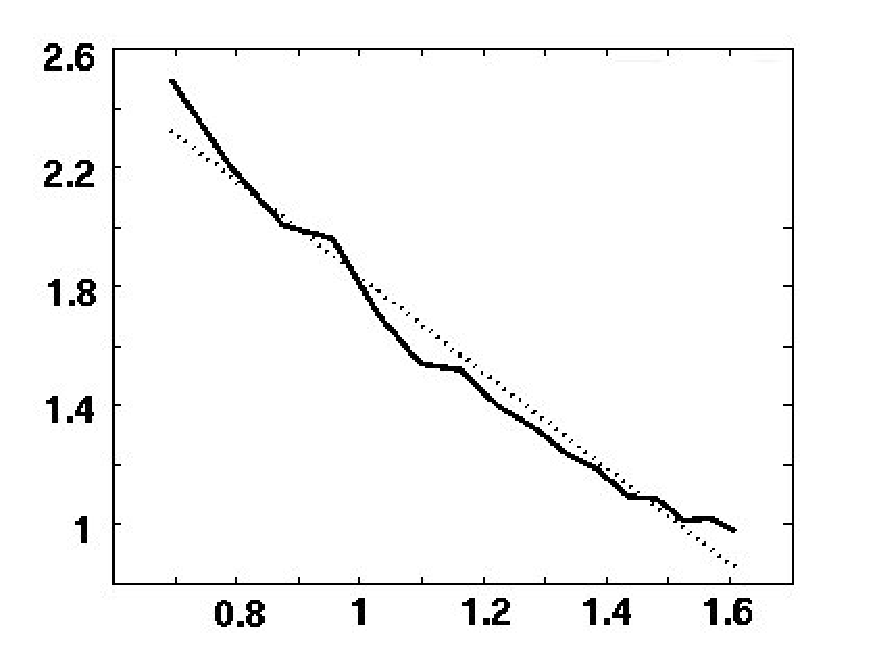}}\\ 
\end{center}
\caption{Density plot of a non-radial non-trapping sound speed (top), density plot of phantoms (left column) and the corresponding $H^1$ errors of reconstruction as functions of the cut-off time $T$ (right column). The plots of the errors are in logarithmic scale. The dotted lines show the linear regression interpolation of the error.  The decay of the error in the first example (second row) is of the order of $T^{-1.68}$, in the second one (third row) it is $T^{-1.61}$.}

\label{F:bump_speed}
\end{figure}
\begin{figure}[!ht]
\begin{center}
\scalebox{0.3}{\includegraphics{trap_crat.eps}}\vspace{0.2cm}\\
\scalebox{0.25}{\includegraphics{vsm_comb1_ph.eps}} \mbox{\hskip 50pt} \scalebox{0.27}{\includegraphics{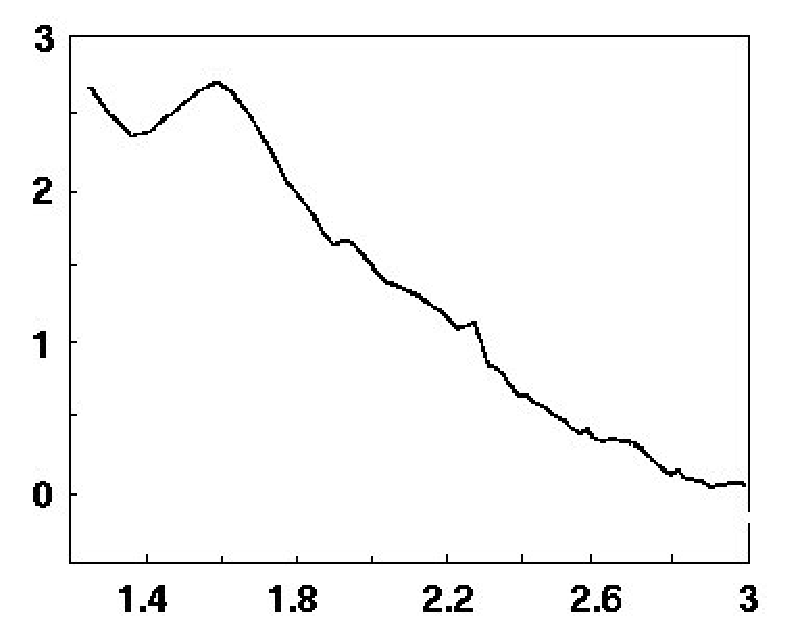}}\vspace{0.2cm}\\
\scalebox{0.25}{\includegraphics{face}} \mbox{\hskip 50pt} \scalebox{0.27}{ \includegraphics{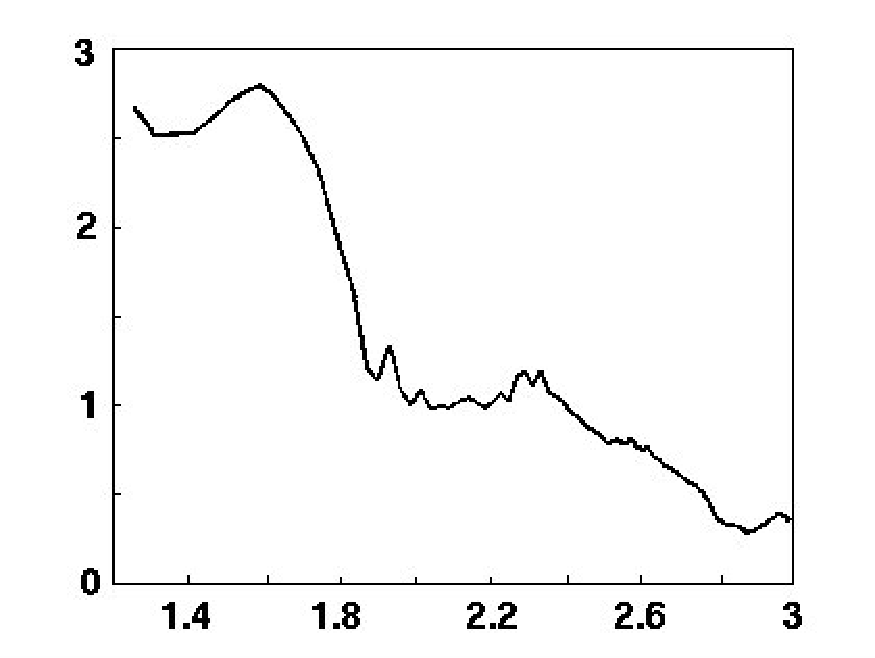}}
\end{center}
\caption{Axial profile of a radial ''trapping crater`` sound speed (top), density plot of phantoms (left column) and the corresponding $L^2$ errors of reconstruction as functions of the cut-off time $T$ (right column). The plots of the errors are in logarithmic scale. }
\label{F:trap_crater}
\end{figure}

\begin{figure}[!ht]
\begin{center}
\scalebox{0.3}{\includegraphics{parab_speed.eps}} \vspace{0.2cm}\\
\scalebox{0.25}{\includegraphics{vsm_comb1_ph.eps}}\mbox{\hskip 50pt} \scalebox{0.27}{\includegraphics{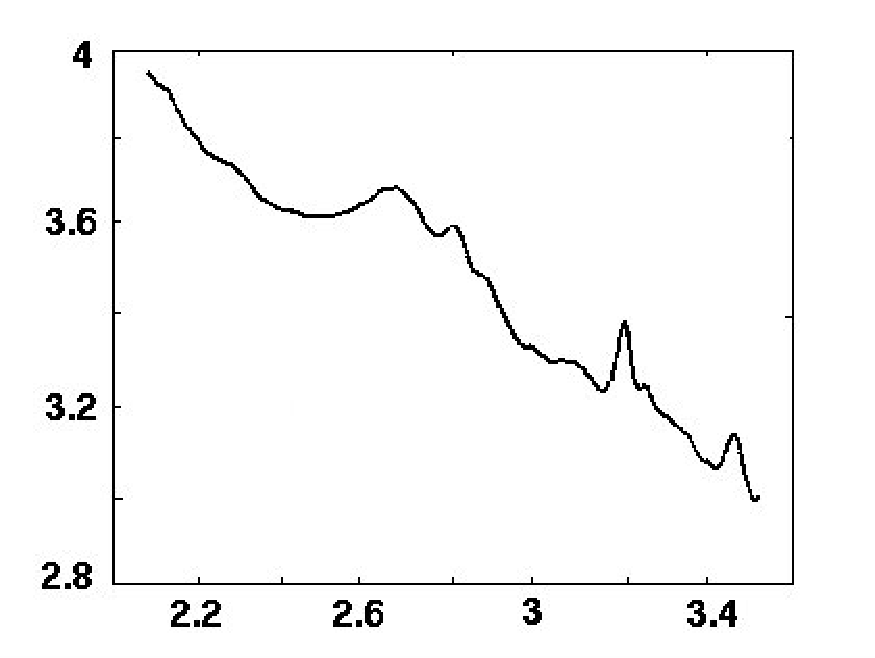}}\\
\end{center}
\caption{Axial profile of a radial parabolic trapping sound speed (top), density plot of a phantom (left column) and the corresponding $L^2$ error of reconstruction as a function of the cut-off time $T$ (right column). The plot of the error is in logarithmic scale.}
\label{F:parab_speed}
\end{figure}
In this section we compare the errors of the numerical reconstructions of several phantoms to their estimates given in Theorem \ref{T:Theorem}. The numerical simulations were done in 2D media with variable sound speed. Both cases of non-trapping and trapping sound speed are considered. 

For the computations we used the rectilinear finite difference scheme (as in \cite{BurgMattHaltPalt,HKN}), implemented in Matlab. For both simulation of the phantom data and reconstruction, we approximated the boundary of the unit circle, $S$, by the set of grid points closest to $S$ and lying within $B$. For the simulation of phantom data the following problem was solved
\[\left\{ \begin{array}{l}
             p_{tt} - c_0^2\Delta p = 0 \qquad \mbox{in} \quad D\times(0,T_0] \\
             p(x,0) = f(x)\;    p_t(x,0)= 0\\
             p|_{\partial D} =0, 
      \end{array}\right.\]
where $D=[-a,a]^2$ was a square containing $S$ and large enough to ensure that no reflections off its boundary would reach $S$ for time $T_0$. The values of the solution on $S$ were recorded for all time steps. The time $T_0$, the domain $D$ and the space and time step-sizes were adjusted depending on the situation. 
For the reconstruction part, the same equation and difference scheme (but a different mesh) were used, this time on the square $[-1.2,1.2]^2$ instead of  on $D$. To obtain the values of the boundary data at the grid points used for reconstruction we used nearest point approximation. In order to observe the behavior of the error of reconstruction with respect to the cut-off time $T$, multiple reconstructions from different cut-off times were made and the error was graphed on a logarithmic scale. In both forward simulation and reconstruction, the spatial step-size $h$ varied case by case and was of the order of $10^{-2}$. This corresponds to using several hundred detectors on the boundary of the unit circle. 

On Figure \ref{F:rec} examples of reconstructions using the time reversal method are shown. In the first case a variable non-trapping speed was used and the reconstruction is of a high quality. In the second example, the radial sound speed equals $|x|$ in the annulus $A=\{x|0.5<|x|<1\}$ and is constant elsewhere. This ''trapping crater'' speed traps all rays corresponding to bicharacteristics that start at points $(x_0,\xi_0)$ such that $x_0\in A$ and $x_0\perp\xi_0$. The trapping leads to the blurring of the radial sides of the phantom, since those singularities do not reach the observation surface. A discussion of this phenomenon can be found in \cite{HKN}. In the third example, shown in Figure \ref{F:rec}, the sound speed equals $|x|^2+0.1$ inside the unit circle. This speed exhibits a more severe trapping than the ``trapping crater`` speed, which leads to a much stronger blurring. In this case, for any $x_0$ there exists a cone of directions $\xi_0$, such that the ray corresponding to the bicharacteristic starting at $x_0,\xi_0$ is trapped \cite{Linh}.

In the next two figures we present various sound speeds and analyse the errors of reconstruction as functions of the cut-off time $T$. These figures show the particular sound speed, phantoms and the corresponding errors, but not actual reconstructions, as we have already seen examples of such in Figure \ref{F:rec}. The plot of each error shows only values of the error after some minimal time, specific to each sound speed, after which it makes sense to apply the time reversal procedure. Namely, this minimal time is taken to be the time a wave needs to cross the unit circle. After the error has decreased to a very small value it levels off, which is due to the discretization of the model. Thus, the error plots do not show these values. In order to easily observe the behaviour of the error, the plots are made in a logarithmic scale and the corresponding linear regression interpolations are graphed.

Figures \ref{F:nice_cos} and \ref{F:bump_speed} show two examples of variable non-trapping sound speeds. In both cases two different phantoms are considered and the $H^1$ error of reconstruction as a function of the cut-off time $T$ is plotted. In Theorem \ref{T:Theorem} we estimated that, in two dimensions, the error behaves as $T^{-1}$. The examples from Figures \ref{F:nice_cos} and \ref{F:bump_speed} suggest even faster decay for those particular sound speeds. In the case shown in Figure \ref{F:nice_cos} the $H^1$ error decays as $T^{-1.3}$. The decay of the error in the second example (Figure \ref{F:bump_speed}) is even faster: the plots suggest  that it is of the order of $T^{-1.6}$. 
The next two examples deal with trapping sound speeds. In these cases the behaviour of the error depends also on the initial condition, as there is no uniform local energy decay unless additional smoothness is assumed (e.g. {\cite{Ralston82,Burq}}). The $L^2$ errors shown in Figures \ref{F:trap_crater} and \ref{F:parab_speed} decay as $T$ increases.

\section{Final Remarks and Conclusion}\label{S:conclusion}

\begin{itemize}
\item In this text, we have concentrated on the errors resulting from the time reversal approximation. One is also interested in the stability of the method with respect to errors in the data. This question is answered by the standard results on stability of the mixed problem for the wave equation. If $n(x,t)$ is the error
in the measured data  $g(x,t)$, then a stability estimate should bound a norm of $\tilde{e}(x,0;T)$, where $\tilde{e}(x,t;T)$ solves the following mixed value problem: 
\begin{equation}\label{E:err_noise}
\left\{ \begin{array}{l}
             \tilde{e}_{tt}(x,t;T) = c^2(x)\Delta \tilde{e}(x,t;T) \qquad \mbox{in} \; B\times[0,T] \\
             \tilde{e}(x,T;T) = 0\\
             \tilde{e}_t(x,T;T)= 0\\
             \tilde{e}|_S(x,t;T)= n_\eg(x,t)\colon=(1-\varphi_\eg(t))n(x,t) \qquad \mbox{in}\;S\times[0,T]. 
         \end{array}\right.
\end{equation}

Since $c(x)$ is assumed to be smooth, classical results (e.g. \cite[Ch. 5, Sec. 5]{LionsII}) give, for instance, the following estimate:
\[
\| \tilde{e}(x,0;T)\|_{L^2(B)} \leq C \|n(x,t)\|_{H^{1+\delta,1+\delta}([0,T]\times S)}
\]
Here $H^{1+\delta,1+\delta}([0,T]\times S)=L^2(0,T;H^{1+\delta}(S))\cap H^{1+\delta}(0,T;L^2(S))$ and $\delta >0$. 
In a nutshell, the reconstruction is as stable as Radon inversion in the corresponding dimension. 
\item Our numerical examples show decay of the error also in the trapping speed case. One could probably get some estimates on this decay, assuming sufficient smoothness of the function $f(x)$ to be reconstructed and using the results of \cite{Burq}. It is unlikely that one can get such estimates without a smoothness assumption.
\end{itemize}
To conclude, in this paper we have provided error estimates of the time reversal justifying it in any dimension under non-trapping condition. The numerical examples agree with the error estimates in the cases of non-trapping speeds and show that the time reversal works even when the speed is trapping. 

\section*{Acknowledgments}
This work was supported by the NSF grants DMS 0715090 and DMS 0604778.
The author is grateful to P. Kuchment for fruitful discussions and advice on this paper, and to L. Nguyen for information. The author also wishes to thank the reviewers for the comments, which helped to improve the paper. 


\end{document}